\documentclass[number,citesort,MSNbibl,secthm,seceqn,dvips]{arxbj}
\usepackage{dcolumn}
\usepackage{graphicx}

\aid{0}
\volume{18}
\issue{2}
\pubyear{2012}
\firstpage{678}
\lastpage{702}
\doi{10.3150/10-BEJ344}

\makeatletter
\newcolumntype{d}[1]{D{.}{.}{#1}}
\newtheorem{lem}{Lemma}[section]
\newproclaim{defn}{Definition}
\newremark{rem}{Remark}[section]
\newremark{exm}{Example}[section]
\newproclaim{ass}{Assumption}[section]

\makeatother

\begin{document}
\begin{frontmatter}

\title{Estimation in semi-parametric regression with non-stationary regressors}
\runtitle{Estimation in semi-parametric non-stationary regression}

\begin{aug}
\author{\fnms{Jia} \snm{Chen}\thanksref{e1}\ead[label=e1,mark]{jia.chen@adelaide.edu.au}},
\author{\fnms{Jiti} \snm{Gao}\corref{}\thanksref{e2}\ead[label=e2,mark]{jiti.gao@adelaide.edu.au}} \and
\author{\fnms{Degui} \snm{Li}\thanksref{e3}\ead[label=e3,mark]{degui.li@adelaide.edu.au}}
\runauthor{J. Chen, J. Gao and D. Li}
\address{School of Economics, The University of Adelaide, Adelaide SA 5005,
Australia.\\ \printead{e1,e2,e3}}
\end{aug}

\received{\smonth{5} \syear{2010}}
\revised{\smonth{10} \syear{2010}}

%
\begin{abstract}
In this paper, we consider a partially linear model of the form
$Y_{t}=X_{t}^\tau\theta_0+g(V_{t})+\epsilon_{t}$, $t=1, \ldots,n$,
where $\{V_{t}\}$ is a $\beta$ null recurrent Markov chain,
$\{X_{t}\}$ is a sequence of either strictly stationary or
non-stationary regressors and $\{\epsilon_{t}\}$ is a stationary
sequence. We propose to estimate both $\theta_0$ and $g(\cdot)$ by
a semi-parametric least-squares (SLS) estimation method.
Under certain conditions, we then show that the proposed SLS
estimator of $\theta_0$ is still asymptotically normal with the
same rate as for the case of stationary time series. In addition, we
also establish
an asymptotic distribution for the nonparametric estimator of the
function $g(\cdot)$. Some numerical examples are provided to show that
our theory and estimation method work well in practice.
\end{abstract}

%

\end{frontmatter}

\section{Introduction}

During the past two decades, there has been much interest in
various nonparametric and semi-parametric techniques to model time
series data with possible nonlinearity. Both estimation and
specification testing problems have been systematically examined
for the case where the observed time series satisfy a type of
stationarity. For more details and recent developments, see
Robinson \cite{Rob83,Rob88,Rob89},
Fan and Gijbels \cite{FanGij96},
H\"{a}rdle \textit{et al.} \cite{Haretal97,HarLiaGao00},
Fan and Yao \cite{FanYao03},
Gao \cite{Gao07},
Li and Racine \cite{LiRac07} and the references therein.

As pointed out in the literature, the stationarity assumption
seems too restrictive in practice. For example, when tackling
economic and financial issues from a time perspective, we often
deal with non-stationary components. In reality, neither prices nor
exchange rates follow a stationary law over time. Thus
practitioners might feel more comfortable avoiding restrictions
like stationarity for processes involved in economic time series
models. There is much literature on parametric linear and
nonlinear models of non-stationary time series, but very little work
has been done in nonparametric and semi-parametric nonlinear cases.
In nonparametric estimation of nonlinear regression and
autoregression of non-stationary time series models and\vadjust{\goodbreak}
continuous-time financial models, existing studies include
Phillips and Park \cite{Phi},
Karlsen and Tj{\o}stheim \cite{KarTjs01},
Bandi and Phillips~\cite{BanPhi03},
Karlsen \textit{et al.} \cite{KarMykTjs07},
Schienle \cite{Sch} and
Wang and Phillips \cite{WanPhi09N1,WanPhi09N2}. Recently,
Gao \textit{et al.}~\cite{Gaoetal09,Gao} considered
nonparametric specification testing in both autoregression and
cointegration models.\looseness=-1

Consider a nonparametric regression model of the form
%
\begin{equation}\label{cgl1.1}
Y_t=m(Z_t)+\epsilon_t,\qquad t=1,\ldots,  n,
\end{equation}
where $\{Y_t\}$ and $\{Z_t\}$ are non-stationary time series,
$m(\cdot)$ is an unknown function defined in $\mathbf{R}^p$ and
$\{\epsilon_t\}$ is a sequence of strictly stationary errors. We
may apply a nonparametric method to estimate $m(\cdot)$,
%
\begin{equation}\label{cgl1.2}
\widehat{m}(z):=\widehat{m}_n(z)=\sum_{t=1}^n a_{nt}(z)
Y_t,
\end{equation}
where $\{a_{nt}(z)\}$ is a sequence of positive weight functions; see
Karlsen \textit{et al.} \cite{KarMykTjs07} and
Wang and Phillips \cite{WanPhi09N1,WanPhi09N2}.

As pointed out in the literature for the case where the dimension
of $\{Z_t\}$ is larger than three, $m(\cdot)$ may not be estimated
by $\widehat{m}(z)$ with reasonable accuracy due to ``the curse of
dimensionality''. The curse of dimensionality problem has been
clearly illustrated in several books, such as
Silverman \cite{Sil86},
Hastie and Tibshirani \cite{HasTib90},
Green and Silverman~\cite{12a},
Fan and Gijbels \cite{FanGij96},
H\"{a}rdle \textit{et al.} \cite{HarLiaGao00},
Fan and Yao \cite{FanYao03} and
Gao \cite{Gao07}. There are several ways to circumvent the curse of
dimensionality. Perhaps one of the most commonly used methods is
semi-parametric modelling, which is taken to mean partially linear
modelling in this context. In this paper, we propose using a
partially linear model of the form
%
\begin{equation}\label{cgl1.3}
Y_{t}=X_{t}^\tau\theta_0+g(V_{t})+\epsilon_{t}, \qquad t=1,\ldots,  n,
\end{equation}
where $\theta_0$ is an unknown $d$-dimensional vector; $g(\cdot)$ is
some continuous function; $\{X_{t}=(x_{t1},\ldots,  x_{td})^\tau\}$
is a sequence of either stationary or non-stationary
regressors, as assumed in A1 below; $\{V_{t}\}$ is a $\beta$ null
recurrent Markov process (see Section~\ref{s2} below for detail); and
$\{\epsilon_{t}\}$ is an error process. As discussed in
Section \ref{s3.2} below,
$\{\epsilon_t\}$ can be relaxed to be either stationary and
heteroscedastic or non-stationary and heteroscedastic.

An advantage of the partially linear approach is that any existing
information concerning possible linearity of some of the
components can be taken into account in such models.
Engle \textit{et al.} \cite{Eng} were among the first to study this kind of partially
linear model. It has been studied extensively in both econometrics
and statistics literature. With respect to development in the
field of semi-parametric time series modelling, various estimation
and testing issues have been discussed for the case where both
$\{X_t\}$ and $\{V_t\}$ are strictly stationary (see, e.g.,
H\"{a}rdle \textit{et al.} \cite{HarLiaGao00} and
Gao \cite{Gao07}) since the publication of
Robinson~\cite{Rob88}. For the case where $\{V_t\}$ is a sequence of
either fixed designs or strictly stationary regressors but there
is some type of unit root structure in $\{X_t\}$, existing
studies, such as Juhl and Xiao \cite{JuhXia05},
have discussed estimation and testing problems.

To the best of our knowledge, the case where either $\{V_t\}$ is a
sequence of non-stationary regressors or both $\{X_t\}$ and
$\{V_t\}$ are non-stationary has not been discussed in the
literature. This paper considers the following two cases: (a)
where $\{X_t\}$ is a sequence of strictly stationary regressors\vadjust{\goodbreak}
and $\{V_t\}$ is a sequence of non-stationary regressors; and (b)~where
both $\{X_t\}$ and $\{V_t\}$ are non-stationary. In this
case, model~(\ref{cgl1.3}) extends some existing models
(Robinson \cite{Rob88},
H\"{a}rdle \textit{et al.} \cite{HarLiaGao00},
Juhl and Xiao \cite{JuhXia05} and
Gao \cite{Gao07})
from the case where $\{V_t\}$ is a sequence of strictly stationary
regressors to the case where $\{V_t\}$ is a sequence of
non-stationary regressors. Since the invariant
distribution of the $\beta$ null recurrent Markov process
$\{V_t\}$ does not have any compact support, however, the semi-parametric
technique used in stationary time series cannot be directly applicable
to our case. In this paper, we will develop a new semi-parametric
estimation method to address such new technicalities when establishing
our asymptotic theory.

The main objective of this paper is to derive asymptotically
consistent estimators for both $\theta_0$ and $g(\cdot)$ involved
in model~(\ref{cgl1.3}). In a traditional stationary time series
regression problem, some sort of stationary mixing condition is
often imposed on the observations $(X_t,V_t)$ to establish
asymptotic theory. In this paper, it is interesting to find that
the proposed semi-parametric least-squares (SLS) estimator of
$\theta_0$ is still
asymptotically normal with the same rate as that in the case of
stationary time series when certain smoothness conditions are
satisfied. In addition,
our nonparametric estimator of $g(\cdot)$ is also asymptotically
consistent, although the rate
of convergence, as expected, is slower than that for the
stationary time series case.

The rest of the paper is organized as follows. The estimation
method of $\theta_0$ and $g(\cdot)$ and some necessary conditions
are given in Section \ref{s2}. The main results and some
extensions are provided in Section \ref{s3}.
Section \ref{s4} provides a simulation
study. An analysis of an economic data set from the United States is
given in Section \ref{s5}. An outline of the proofs of the main theorems is
given in Section \ref{s6}.
Supplementary Material section
gives a~description for a supplemental
document by Chen, Gao and Li \cite{CheGaoLi10}, from which the detailed proofs of
the main theorems, along with some technical lemmas, are available.

\section{Estimation method and assumptions}\label{s2}

\subsection{Markov theory}\label{s2.1}

Let $\{V_{t},t\geq0\}$ be a Markov chain with transition
probability $P$ and state space $(\mathbf{E},\mathcal{E})$, and
$\phi$ be a measure on $(\mathbf{E},\mathcal{E})$. Throughout the
paper, $\{V_{t}\}$ is assumed to be $\phi$-irreducible Harris
recurrent, which makes asymptotics for semi-parametric estimation
possible. The class of stochastic
processes we are dealing with in this paper is not the general
class of null recurrent Markov chains. Instead, we need to impose
some restrictions on the tail behavior of the distribution of the
recurrence time $S_{\alpha}$ of the chain. This is what we are
interested in: a class of
$\beta$ null recurrent Markov chains.
\begin{defn*}
A Markov chain $\{V_{t}\}$ is
$\beta$ null recurrent if there exist a small non-negative
function $f(\cdot)$ (the definition of a small function can be found in
the supplemental document), an initial measure $\lambda$, a constant
$\beta\in(0,1)$ and a slowly varying function $L_{f}(\cdot)$ such that
%
\begin{equation}\label{cgl2.1}
\mathbf{E}_{\lambda}\Biggl[\sum_{t=0}^{n}f(V_{t}) \Biggr]
\sim \frac{1}{\Gamma(1+\beta)} n^{\beta}L_{f}(n) \qquad
\mbox{as }n\rightarrow\infty,\vadjust{\goodbreak}
\end{equation}
where $\mathbf{E}_{\lambda}$ stands for the expectation with initial
distribution $\lambda$ and $\Gamma(\cdot)$ is the usual gamma
function.\vspace*{-3pt}
\end{defn*}

It is shown in Karlsen and Tj\o{}stheim \cite{KarTjs01} that when there
exist some small measure~$\nu$ and small function $s$ with
$\nu(\mathbf{E})=1$ and $0\leq s(v)\leq1$, $v\in\mathbf{E}$, such that
%
\begin{equation}\label{cgl2.2}
P\geq s\otimes\nu,\vspace*{-1pt}
\end{equation}
then $\{V_{t}\}$ is $\beta$ null recurrent if and only if
%
\begin{equation}\label{cgl2.3}
\mathbf{P}_{\alpha}(S_{\alpha}>n)=\frac{1}{\Gamma(1-\beta)n^{\beta}L_{s}(n)}
\bigl(1+\mathrm{o}(1)\bigr),\vspace*{-1pt}
\end{equation}
where $L_{s}=\frac{L_{f}}{\pi_{s}f}$ and $\pi_s$ is the invariant
measure as defined in Karlsen and Tj{\o}stheim \cite{KarTjs01}.
Furthermore, if
(\ref{cgl2.3}) holds, by Lemma 3.4 in Karlsen and Tj{\o}stheim \cite{KarTjs01},
$\widehat{\beta}:=\frac{\ln(N_{\mathbf{C}}(n))}{\ln n}$ is a strongly
consistent estimator of $\beta$, where $N_{\mathbf{C}}(n)=\sum
_{t=1}^n I_{\mathbf{C}}(V_t)$, in which $I_A(\cdot)$ is the
conventional indicator function and ${\mathbf{C}}$ is a small set as
defined in Karlsen and Tj{\o}stheim~\cite{KarTjs01}.

We then introduce a useful decomposition that is critical in the
proofs of asymptotics for nonparametric estimation in null
recurrent time series. Let $f$ be a real function defined in $\mathbf{R}$. We now decompose the partial sum
$S_{n}(f)=\sum_{t=0}^{n}f(V_{t})$ into a sum of independent and
identically
distributed (i.i.d.) random variables with one main part and two
asymptotically negligible minor parts. Define {\small
\[
Z_{k}= \cases{%
\displaystyle \sum_{t=0}^{\tau_{0}}f(V_{t}),&\quad $k=0$,\vspace*{-1pt}\cr
\displaystyle\sum_{t=\tau_{k-1}+1}^{\tau_{k}}f(V_{t}),&\quad $1\leq k\leq N(n)$,\vspace*{-1pt}\cr
\displaystyle\sum_{t=\tau_{N(n)}+1}^{n}f(V_{t}),&\quad $k=(n)$,}\vspace*{-1pt}
\]
where} the definitions of $\tau_k$ and $N(n)$ will be given in the
supplemental document. Then
%
\begin{equation}\label{cgl2.4}
S_{n}(f)=Z_{0}+\sum_{k=1}^{N(n)}Z_{k}+Z_{(n)}.\vspace*{-1pt}
\end{equation}
From Nummelin's \cite{22a}
result, we know that $\{Z_{k}, k\geq1\}$ is
a sequence of i.i.d. random variables. In the decomposition
(\ref{cgl2.4}) of $S_{n}(f)$, $N(n)$ plays the role of the number
of observations. It follows from Lemma 3.2 in
Karlsen and Tj\o stheim \cite{KarTjs01}
that $Z_{0}$ and $Z_{(n)}$ converge to zero almost
surely when they are divided by $N(n)$. Furthermore,
Karlsen and Tj\o{}stheim \cite{KarTjs01}
show that if (\ref{cgl2.2}) holds and $\int
|f(v)|\pi_{s}(\mathrm{d}v)<\infty$, then for an arbitrary initial
distribution $\lambda$ we have
%
\begin{equation}\label{cgl2.5}
\frac{1}{N(n)}S_{n}(f)\longrightarrow\pi_{s}(f)
\qquad \mbox{almost surely (a.s.)},\vspace*{-1pt}
\end{equation}
where $\pi_{s}(f)=\int f(v)\pi_{s}(\mathrm{d}v)$.\vadjust{\goodbreak}

Some useful results for Markov theory are available from Appendix A of
the supplemental document.\vspace*{-3pt}

\subsection{Estimation method}\label{s2.2}\vspace*{-3pt}

As assumed in assumption A1 below, there exist a function $H(\cdot
)$ and a stationary process~$\{U_t\}$ such that $X_t = H(V_t) + U_t$.
Since $E[\epsilon_t|V_t=v]=E[\epsilon_t]=0$ is assumed in A2(ii)
and~A3(ii), we have
%
\begin{equation}\label{cgl2.6a}
E[Y_t|V_t=v] = E \bigl[ \bigl(X_t^{\tau} \theta_0 + g(V_t) + \epsilon_t \bigr)|V_t=v \bigr]
= H(v)^{\tau}\theta_0 + g(v).
\end{equation}

This implies that $\Psi(v) \equiv E[Y_t|V_t=v]$ is a function of $v$
independent of $t$ for each fixed~$v$ and given $\theta_0$. Thus, the
form of $g(v)$ can be represented by
%
\begin{equation}\label{cgl2.6}
g(v) = \Psi(v) - H(v)^{\tau} \theta_0.
\end{equation}

In view of (\ref{cgl2.6}), we can rewrite model (\ref{cgl1.3}) as
%
\begin{equation}\label{cgl2.7}
Y_t - \Psi(V_t) =  \bigl(X_t - H(V_t) \bigr)^{\tau} \theta_0 + \epsilon_t.
\end{equation}

Letting $W_t = Y_t - \Psi(V_t)$ and $U_t = X_t - H(V_t)$,
model (\ref{cgl2.7}) implies
%
\begin{equation}\label{cgl2.9}
W_t =Y_t - \Psi(V_t) =  \bigl(X_t - H(V_t) \bigr)^{\tau} \theta
_0+\epsilon_t = U_t^{\tau} \theta_0 + \epsilon_t.
\end{equation}

Note that $E[W_t]=E[U_t^{\tau}\theta_0]+E[\epsilon_t]=0$. In the case
where $\{(X_t, V_t, \epsilon_t)\dvt t\geq1\}$ is a~sequence of
stationary random variables, various estimation methods for $\theta_0$
and $g(\cdot)$ in model~(\ref{cgl1.3}) have been
studied by many authors (see, e.g., Robinson~\cite{Rob88},
H\"{a}rdle \textit{et al.}~\cite{HarLiaGao00}
and Gao~\cite{Gao07}).

We now propose an SLS estimation method
based on the kernel smoothing. For every given $\theta$, we define
a kernel estimator of $g(v)$ by
%
\begin{equation}\label{eq2.1}
g_{n}(v;\theta)=\sum_{t=1}^{n} w_{nt}(v)(Y_{t}-X_{t}^\tau\theta),
\end{equation}
where $\{w_{nt}(v)\}$ is a sequence of weight functions given by
\[
w_{nt}(v)=\frac{K_{v,h}(V_{t})}{\sum_{k=1}^{n} K_{v,h}(V_{k})}
\qquad  \mbox{with }
K_{v,h}(V_t) = \frac{1}{h} K \biggl(\frac{V_t-v}{h} \biggr),
\]
in which $K(\cdot)$ is a probability kernel function and $h=h_n$ is a
bandwidth parameter.

Replacing $g(V_{t})$ by $g_{n}(V_{t};\theta)$ in model (\ref{cgl1.3})
and applying the SLS estimation method, we obtain the SLS estimator,
$\overline{\theta}_n$, of $\theta_0$ by minimizing
\[
\frac{1}{n} \sum_{t=1}^n
\bigl(Y_t - X_t^{\tau} \theta- g_n(V_t, \theta) \bigr)^2
\]
over $\theta$. This implies
%
\begin{equation}\label{cgl2.10}
\overline{\theta}_n=(\overline{X}^\tau\overline{X})^{-1}
\overline{X}^\tau\overline{Y},\vadjust{\goodbreak}
\end{equation}
where
$\overline{X}^\tau=(\widetilde{X}_{1},\ldots,\widetilde{X}_{n})$,
$\widetilde{X}_{t}=X_{t}-\sum_{k=1}^{n}w_{nk}(V_{t})X_{k}$,
$\overline{Y}^\tau=(\widetilde{Y}_{1},\ldots,\widetilde{Y}_{n})$
and $\widetilde{Y}_{t}=Y_{t}-\sum_{k=1}^{n}
w_{nk}(V_{t})Y_{k}$. And $g(\cdot)$ is then estimated by
%
\begin{equation}\label{cgl2.11}
\overline{g}_n(\cdot)=g_{n}(\cdot;\overline{\theta}_n).
\end{equation}

This kind of estimation method has been studied in the literature
(see, e.g., H\"{a}rdle \textit{et al.} \cite{HarLiaGao00}).
When $\{V_t\}$ is
a sequence of either fixed designs or stationary regressors with a
compact support, the conventional weighted least-squares estimators
(\ref{cgl2.10}) and (\ref{cgl2.11}) work well in both the large and
small sample cases.
Since the invariant distribution of~$\beta$ null recurrent Markov chain
$\{V_t\}$
might not have any compact support, it is difficult to
establish asymptotic results for the estimators (\ref{cgl2.10})
and (\ref{cgl2.11}) owing to the random denominator problem
involved in $w_{nt}(\cdot)$. Hence, to establish
our asymptotic theory, we apply the following weighted least-squares
estimation method (see, e.g., Robinson \cite{Rob88}).
Define
%
\begin{equation}\label{cgl2.12}
F_t:=F_{nt}=I \bigl(|p_n(V_t)|>b_n \bigr),
\end{equation}
where
\[
p_n(v)=\frac{1}{N(n)}\sum_{k=1}^n K_{v,h}(V_k)
\]
and $\{b_n\}$ is a sequence of positive numbers satisfying some
conditions. Furthermore, let
\[
\widetilde{X}^\tau=(\widetilde{X}_{1}F_1,\ldots,\widetilde{X}_{n}F_n)
\quad \mbox{and} \quad
\widetilde{Y}^\tau=(\widetilde{Y}_{1}F_1,\ldots,\widetilde{Y}_{n}F_n).
\]
Throughout this paper, we propose to estimate $\theta_0$ by
%
\begin{equation}\label{cgl2.13}
\widehat{\theta}_n=(\widetilde{X}^\tau\widetilde{X})^{-1}
\widetilde{X}^\tau\widetilde{Y}
\end{equation}
and $g(\cdot)$ by
%
\begin{equation}\label{cgl2.14}
\widehat{g}_n(\cdot)=g_{n}(\cdot;\widehat{\theta}_n).
\end{equation}

\subsection{Assumptions}\label{s2.3}

As may be seen from equation (\ref{cgl2.9}), further discussion on
the semi-parametric estimation method depends heavily on the
structure of $\{X_t\}$ and $\{V_t\}$. This paper is concerned with the
following two cases: (i) where
$\{X_t\}$ is a sequence of strictly stationary regressors and
independent of $\{V_t\}$; and (ii) where $\{X_t\}$ is a sequence
of non-stationary regressors with the non-stationarity being generated
by $\{V_t\}$.

Before stating the main assumptions, we introduce the definition
of $\alpha$ mixing dependence. The stationary sequence $\{Z_t,
t=0,\pm1,\ldots\}$ is said to be $\alpha$ mixing if
$\alpha(n)\rightarrow0$ as $n\rightarrow\infty$, where
\[
\alpha(n)=\sup_{A\in{\mathcal{F}}_{-\infty}^{0},
B\in{\mathcal{F}}_{n}^{\infty}}|P(AB)-P(A)P(B)|,
\]
in which $\{{\mathcal{F}}_k^j\}$ denotes a sequence of
$\sigma$ fields generated by $\{Z_t,k\leq t\leq j\}$. Since its
introduction by Rosenblatt \cite{Ros56}, $\alpha$ mixing dependence is
a property shared by many time series models (see, e.g.,
Withers \cite{Wit81} and
Gao \cite{Gao07}). For more details
about limit theorems for $\alpha$ mixing processes, we refer to
Lin and Lu \cite{LinLu96} and the references therein.

The following assumptions are necessary to derive the asymptotic
properties of the semi-parametric estimators.
\begin{longlist}[(A4.)]
\item[A1.] There exist an unknown function $H(v)$ and a stationary
process $\{U_t\}$ such that $X_t = H(V_t) + U_t$.
\item[A2.] (i) Suppose that $\{U_t\}$ is a stationary ergodic Markov process
with $E[U_1]=0$ and $E [\Vert
U_1\Vert^{4+\gamma_1} ]<\infty$ for some $\gamma_1>0$, where
$\|\cdot\|$ stands for the Euclidean norm. Furthermore, we suppose
that $\Sigma:=E[U_1U_1^\tau]$ is positive definite and $\{U_t\}$
is $\alpha$ mixing with
%
\begin{equation}\label{cgl2.16}
\sum_{t=1}^\infty\alpha_U^{{\gamma_1}/{(4+\gamma_1)}}(t)<\infty,
\end{equation}
where $\alpha_U(t)$ is the $\alpha$ mixing coefficient of
$\{U_t\}$.

(ii) Let $\{\epsilon_t\}$ be a stationary ergodic Markov process
with $E[\epsilon_1]=0$, $\sigma^2:=E [\epsilon_1^2 ]>0$
and $E [|\epsilon_1|^{2+\gamma_2} ]<\infty$ for some
$\gamma_2>0$. Furthermore, the process $\{\epsilon_t\}$ is
$\alpha$ mixing with
%
\begin{equation}\label{cgl2.17}
\sum_{t=1}^\infty\alpha_{\epsilon}^{{\gamma_2}/{(2+\gamma_2)}}(t)<\infty,
\end{equation}
where $\alpha_\epsilon(t)$ is the $\alpha$ mixing coefficient of $\{
e_t\}$.
\item[A3.] (i) The invariant measure $\pi_{s}$ of the
$\beta$ null recurrent Markov chain $\{V_t\}$ has a~uniformly
continuous density function $p_{s}(\cdot)$.

(ii) Let $\{U_t\}$, $\{V_t\}$ and $\{\epsilon_t\}$ be mutually
independent.
\item[A4.] Let $f_{i,k}(\cdot)$ be the density function of
\[
V_{i,k}=\varphi_{i-k} (V_{i}-V_{k} ) \qquad \mbox{for $i>k$ with }
\varphi_m=m^{\beta-1}L_s(m) \mbox{ for }m\geq1.
\]
Let
%
\begin{equation}\label{cgl2.18}
\inf_{\delta>0}\limsup_{m\rightarrow\infty}
\sup_{i\geq1}\sup_{|v|\leq\delta}f_{i+m,i}(v)<\infty.
\end{equation}
Furthermore, there exists a sequence of $\sigma$ fields
$\{\mathcal{F}_{t},t\geq0\}$ such that $\{V_{t}\}$ is adapted to~$\mathcal{F}_{t}$.
With probability 1,
%
\begin{equation}\label{cgl2.19}
\inf_{\delta>0}\limsup_{m\rightarrow\infty}
\sup_{i\geq1}\sup_{|v|\leq\delta}f_{i+m,i}(v|\mathcal{F}_{i})<\infty,
\end{equation}
where $f_{i,k}(v|{\mathcal{F}_{k}})$ is the conditional density
function of $V_{i,k}$ given $\mathcal{F}_{k}$.
\item[A5.] (i) The function $g(v)$ is differentiable and the derivative
is continuous in $v\in\mathbf{R}$. In addition, for $n$ large enough
%
\begin{equation}\label{cgl2.20}
\sum_{t=1}^n\int (g^{\prime}(\varphi_t^{-1}v)
)^2f_{t,0}(v)\,\mathrm{d}v=\mathrm{O}(nh^{-1}),
\end{equation}
where $g^{\prime}(\cdot)$ is the derivative of $g(\cdot)$, the
definitions of $\varphi_t$ and $f_{t,0}(v)$ are given in A4
above.\looseness=1

(ii) The function $H(v)$ is differentiable and the derivative is
also continuous in $v\in\mathbf{R}$. In addition, for $n$ large
enough
%
\begin{eqnarray}\label{cgl3.3}
\sum_{t=1}^n\int\|H^{\prime}(\varphi_t^{-1}v)\|^2
f_{t,0}(v)\,\mathrm{d}v &=& \mathrm{O}(nh^{-1}) 
\end{eqnarray}
and
\begin{eqnarray}\label{cgl3.4}
\sum_{t=1}^n\int \|g^{\prime}(\varphi_t^{-1}v)
H^{\prime}(\varphi_t^{-1}v) \|f_{t,0}(v)\,\mathrm{d}v
&=& \mathrm{O} (n^{1/2-\varepsilon_1}b_n^{2}h^{-2} ),
\end{eqnarray}
where $\varepsilon_1>0$ is small enough.
\item[A6.] (i) The probability kernel function $K(\cdot)$ is a
continuous and
symmetric function having some compact support.

(ii) The sequences $\{h_n\}$ and $\{b_n\}$ both satisfy
as $n\rightarrow\infty$
%
\begin{equation}\label{cgl2.21}
h_n \rightarrow0, \qquad
b_n\rightarrow0,\qquad
n^{\varepsilon_0}h_n b_n^{-4}\rightarrow0
\quad \mbox{and}\quad
n^{\beta-\varepsilon_0}h_n b_n^4\rightarrow\infty
\end{equation}
for some $0<\varepsilon_0<\frac{\beta}{2}$. Moreover,
%
\begin{equation}\label{cgl2.22}
\sum_{t=1}^nP \bigl(p_n(V_t)\leq b_n \bigr)=\mathrm{o}(n).
\end{equation}
\end{longlist}
\begin{rem}\label{rem2.1}
(i) While some parts of assumptions A1--A3 may
be non-standard, they are justifiable in many situations.
Condition A1 assumes that $\{X_t\}$ is generated by $X_t = H(V_t)
+ U_t$. This is satisfied when the conditional mean function $H(v) =
E[X_t|V_t=v]$ exists. In this case, A1 holds automatically with
$U_t = X_t - E[X_t|V_t]$. Condition A1 is also commonly used in
the stationary case (see, e.g., Linton \cite{Lin95}). There are
various examples in this kind of situation (see, e.g., in the
univariate case where $X_t = V_t + \varepsilon_t$, in which $\{
\varepsilon_t\}$ is a sequence of i.i.d. errors with $E[\varepsilon
_t]=0$ and $E[\varepsilon_t^2]<\infty$, and independent of $\{V_t\}$.
In this case, $H(v) = E[X_t|V_t=v] = v$ and $U_t = \varepsilon_t$). As
a consequence, condition~A1 does not include the case where $\{
X_t\}$ is a random walk sequence of the form $X_t = X_{t-1} + \zeta_t$.
Note that the case where the non-stationarity in both $\{X_t\}$ and $\{
V_t\}$ is generated by a common random walk structure will need to be
discussed separately, since the methodology involved is likely to be
quite different. In Section~\ref{s3.2} below, we will give some discussion
about the case where $\{H(V_t)\}$ is replaced by a bivariate function
of the form $\{H(V_t, t)\}$ to take into account the inhomogeneous case.

(ii) The stationarity assumption on $\{U_t\}$ is to ensure that the
conventional $\sqrt{n}$-rate of convergence is achievable and thus it
is possible to construct an asymptotically efficient estimator for
$\theta_0$. The stationarity condition on $\{U_t\}$ also requires that
$X_t$ can be decomposed into a non-stationary component represented by
$H(V_t)$ and a stationary component $\{U_t\}$. The $\alpha$ mixing
dependence in~A2 is a mild condition on $\{U_t\}$ and the errors
process $\{e_t\}$. Karlsen \textit{et al.} \cite{KarMykTjs07} have made similar
assumptions. As discussed in Section~\ref{s3.2} below, A2(i) can be
relaxed to allow for the inclusion of both endogeneity and
heteroscedasticity. Note that A2(ii) can also be relaxed to allow
for the inclusion of a deterministic function in model (\ref{cgl1.3}).
In such cases, model (\ref{cgl1.3}) can be naturally extended to a
semi-parametrc additive model of the form $Y_t = X_t^{\tau} \theta_0 +
g(V_t) + \lambda(U_t, t) + e_t$ as discussed in Section \ref{s3.2} below.

(iii) As we can see from the asymptotic theory below, the condition on
the existence of the inverse matrix $\Sigma^{-1}$ is required in
Theorem~\ref{thm3.1}. In the case where $\{(X_t, V_t)\}$ is a vector of either
independent regressors or stationary time series regressors, H\"{a}rdle
\textit{et al.}~\cite{HarLiaGao00}
also assume similar conditions (see Section 1.3 in
their book) for establishing the asymptotic results for the
conventional least-squares estimators of $\theta_0$ in (\ref{cgl2.10})
and of $g(\cdot)$ in (\ref{cgl2.11}). Condition A3(i) corresponds
to analogous conditions on the density function in the stationary case.
A3(ii) imposes the mutual
independence to avoid involving some extremely technical conditions.
\end{rem}
\begin{rem}\label{rem2.2}
A4 is similar to but weaker than Assumption
2.3(ii) in Wang and~Phil\-lips~\cite{WanPhi09N1}. It is easy to
check that (\ref{cgl2.18}) and (\ref{cgl2.19}) are satisfied with
$\beta=1$ and $L_s(\cdot)\equiv1$ when~$\{V_t\}$ is a sequence of
either i.i.d. or stationary dependent variables. Consider the
random walk case defined by
%
\begin{equation}\label{cgl2.23}
V_t=V_{t-1}+v_t, \qquad t=1,2,\ldots,   V_0=0,
\end{equation}
where $\{v_t\}$ is a sequence of i.i.d. random variables. The
random walk model (\ref{cgl2.23}) is very important in economics
and finance and has been studied by many authors. It corresponds
to a $1/2$ null recurrent process and it is easy to check that
(\ref{cgl2.18}) and (\ref{cgl2.19}) are satisfied with
$\beta=1/2$, $L_s(n)\equiv1$ and ${\mathcal{F}}_k=\sigma(v_i,i\leq
k)$. On the other hand, (\ref{cgl2.18}) and (\ref{cgl2.19}) can be
formulated in terms of the transition probability. For example,
assume that the transition probability of the Markov process
$\{V_t\}$ is defined by
\[
P(x,\mathrm{d}y)=f(x|y)\,\mathrm{d}y.
\]
Let $f^{k}(\cdot)$ be the marginal density of $\{V_k\}$ and
$f^m(x|y)$ be the $m$ step transition density. Then
\[
f_{i+m,i}(v)=\varphi_m^{-1}\int f^m(\varphi_m^{-1}v+y|y)f^i(y)\,\mathrm{d}y,
\]
where $\varphi_m$ is defined in A4.
\end{rem}
\begin{rem}\label{rem2.3}
(i) A5(i) is assumed to make sure that the
bias term of the nonparametric estimator is negligible when
establishing the asymptotic distribution of the semi-parametric
estimator $\widehat{\theta}_n$. When $\{V_t\}$ is the random walk
process defined by (\ref{cgl2.23}), condition A5(i) can be
verified. If
%
\begin{equation}\label{cgl2.24}
g(v)=\varrho_0+\varrho_1v+\varrho_2|v|^{1+\delta_0},
\qquad  0<\delta_0<1/2,
\end{equation}
$n^{\delta_0}h=\mathrm{O}(1)$ and
$f_{t,0}(v)=\mathrm{O}(v^{-(1+2\delta_0+\varsigma)})$ for some
$\varsigma>0$ as $t\rightarrow\infty$ and $v\rightarrow\infty$, we
can show that by A4,
\[
\sum_{t=1}^n\int (g^{\prime}(\varphi_t^{-1}v))^2f_{t,0}(v)\,\mathrm{d}v
=\mathrm{O} \Biggl(\sum_{t=1}^n\varphi_t^{-2\delta_0} \Biggr)
=\mathrm{O}(n^{1+\delta_0}),
\]
which implies (\ref{cgl2.20}).

(ii) Similarly, condition A5(ii) is also verifiable. Consider the
case where
\[
g(v)=\varrho_0+\varrho_1v    \quad \mbox{and} \quad
H(v)=\mathbf{a}_0+\mathbf{a}_1v+\mathbf{a}_2|v|^{1+\delta_1},
\qquad  1<\delta_1<1/2,
\]
in which $\mathbf{a}_k$, $k=0,1,2$, are $d$-dimensional vectors,
$n^{{1}/{2}+\delta_1-\varepsilon_1}h^2=\mathrm{O}(1)$
($\varepsilon_1<\frac{1}{2}-\delta_1$) and
$f_{t,0}(v)=\mathrm{O}(v^{-(1+2\delta_1+\varsigma)})$ for some
$\varsigma>0$ as $t\rightarrow\infty$ and $v\rightarrow\infty$. We
can also show that~(\ref{cgl3.3}) and~(\ref{cgl3.4}) hold for the
random walk case. The detailed calculation is similar to that in
Remark~\ref{rem2.3}(i) above.
\end{rem}
\begin{rem}\label{rem2.4}
(i) Condition A6(i) is a quite natural
condition on the kernel function and has been used by many authors for
the stationary
time series case. The first part of~A6(i) requires that the rate
of $b_n^{-4}\rightarrow\infty$ is slower than that of $n^{\varepsilon
_0} h\rightarrow0$ and the rate of $b_n^{4}\rightarrow0$ is slower
than that of $n^{\beta-\varepsilon_0} h\rightarrow\infty$. Such
conditions are satisfied in various cases. Letting
$b_n = c_b \log^{-1}(n)$ and $h_n=c_h n^{-\zeta_0}$
for some $c_b>0$, $c_h>0$ and
$\varepsilon_0<\zeta_0<\beta- \varepsilon_0$, then the first part of
A6(ii) holds automatically.

(ii) The second part of A6(ii) is imposed to ensure that the
truncated procedure works in this kind of problem. When $\{V_t\}$ is a
sequence of i.i.d. random variables having some compact support~$S$, it
is easy to show that (\ref{cgl2.22}) holds if
$\inf_{x\in S}p(x)>0$, where $p(\cdot)$ is the density function of $\{V_t\}$. In
the case where $\{V_t\}$ is an i.i.d. sequence without any compact
support, Robinson \cite{Rob88} gives
different conditions such that (\ref{cgl2.22}) holds. We can show that
condition A6(ii) is verifiable when $\{V_t\}$ is a random walk
model of the form~(\ref{cgl2.23}). Since the verification is quite
technical, the details are given in the last part of Appendix C in the
supplemental document.
\end{rem}

\section{The main results and their extensions}\label{s3}

\subsection{Asymptotic theory}\label{s3.1}

We now establish an asymptotic distribution of the estimate
$\widehat{\theta}_n$ in the following theorem. The following theorem
includes two cases: (a)
$\{V_t\}$ is a sequence of non-stationary regressors and $\{X_t\}$
is a sequence of strictly stationary regressors and is independent
of~$\{V_t\}$; and (b) both $\{X_t\}$ and $\{V_t\}$ are non-stationary.
\begin{thm}\label{thm3.1}
Let \textup{A1--A5(i)} and \textup{A6} hold. In addition, suppose that
$\Sigma_{\epsilon,U}:=\sigma^2\Sigma+ 2\sum_{t=2}^\infty
E[\epsilon_1\epsilon_t] E[U_1U_t^\tau]$ is positive definite.
\begin{longlist}[(ii)]
\item[(i)] If $\{X_t\}$ is strictly stationary and independent of
$\{V_t\}$, then as $n\rightarrow\infty$,
%
\begin{equation}\label{cgl3.2}
\sqrt{n}(\widehat{\theta}_n-\theta_0)\stackrel{d}{\longrightarrow}
N (0,\Sigma^{-1}\Sigma_{\epsilon,U}\Sigma^{-1} ).
\end{equation}
\item[(ii)] Suppose that both $\{X_t\}$ and $\{V_t\}$ are non-stationary. If,
in addition, \textup{A5(ii)} is satisfied, then (\ref{cgl3.2}) still holds.
\end{longlist}
\end{thm}
\begin{rem}\label{rem3.1}
(i) Theorem \ref{thm3.1} shows that the standard normality
can still be an asymptotic distribution of the SLS estimate even when
non-stationarity is involved. Theorem~\ref{thm3.1}(ii) further shows that the
conventional rate of $\sqrt{n}$ is still achievable when the
non-stationarity in $\{X_t\}$ is purely generated by $\{V_t\}$ and
certain conditions are imposed on the functional forms of $H(\cdot)$
and $g(\cdot)$.

(ii) Since the asymptotic distribution and
asymptotic variance in (\ref{cgl3.2}) are mainly determined by the
stationary sequences $\{\epsilon_t\}$ and $\{U_t\}$, the above
conclusion extends Theorem 2.1.1 of H\"ardle \textit{et al.} \cite{HarLiaGao00}
for the case when $\{X_t\}$, $\{V_t\}$ and $\{\epsilon_t\}$ are
all strictly stationary. In addition, when $\{X_t\}$ is assumed to be strictly
stationary and independent of $\{V_t\}$ in Theorem \ref{thm3.1}(i), the
covariance matrix reduces to the covariance matrix of~$\{X_t\}$ of the
form $\Sigma= E [ (X_1 - E[X_1] )  (X_1 - E[X_1])^{\tau} ]$.
\end{rem}
\begin{rem}\label{rem3.2}
(i) Theorem \ref{rem3.1} establishes an asymptotically
normal estimator for $\theta_0$. As in the independent and stationary
sample case, an interesting issue is how to construct an asymptotically
efficient estimator for $\theta_0$. As discussed in
Chen \cite{Che88} and
H\"{a}rdle \textit{et al.}~\cite{HarLiaGao00},
it can be shown that $\widehat{\theta}$
achieves the smallest possible variance of $\sigma^2 \Sigma^{-1}$ when
both $\{U_t\}$ and $\{\epsilon_t\}$ are independent and $\epsilon_t
\sim N(0, \sigma^2)$.

(ii) Since the publication of the book by
Bickel \textit{et al.} \cite{Bicetal93},
there has been an increasing interest in the field of asymptotic
efficiency in semi-parametric models. There are certain types of
asymptotic efficiency in this kind of semi-parametric setting. H\"
{a}rdle \textit{et al.} \cite{HarLiaGao00}
consider several types of asymptotically
efficient estimators in Chapters 2 and~5 of the book.
Linton \cite{Lin95}
considers second-order efficiency. Bhattacharya and Zhao \cite{BhaZha97}
establish an asymptotically efficient estimator without requiring
finite variance. Chen \cite{Che07} discusses asymptotic efficiency in
nonparametric and semi-parametric models using sieve estimation.

(iii) As shown in the literature, the establishment of an
asymptotically efficient estimator in this kind of semi-parametric
setting requires the availability of uniform convergence of
nonparametric estimation. Since such uniform convergence results are
not readily available and applicable in this kind of non-stationary
situation, we wish to establish some necessary uniform convergence
results first before we may be able to address the issue of asymptotic
efficiency in future research.
\end{rem}

An asymptotic distribution of $\widehat{g}_n(x)$ is given in Theorem
\ref{thm3.2} below.
\begin{thm}\label{thm3.2}
\textup{(i)} Let the conditions of Theorem \ref{thm3.1}\textup{(i)} hold. If, in addition,
$g(\cdot)$ is twice differentiable and the second derivative,
$g^{\prime\prime}(v)$, is
continuous in $v$ and $n^{\beta/5+\varepsilon} h=\mathrm{o}(1)$ for some
$\varepsilon>0$, then as $n\rightarrow\infty,$
%
\begin{equation}\label{cgl3.5}
\sqrt{\sum_{t=1}^n K \biggl(\frac{V_t - v}{h} \biggr)}
\bigl (\widehat{g}_n(v)-g(v) \bigr)
\stackrel{d}{\longrightarrow}N
\biggl(0,\sigma^{2}\int K^2(u)\,\mathrm{d}u \biggr).
\end{equation}

\textup{(ii)} Let the conditions of Theorem \ref{thm3.1}\textup{(ii)}
hold. If, in addition,
$g(\cdot)$ is twice differentiable and the second
derivative, $g^{\prime\prime}(v)$, is continuous in $v$ and
$n^{\beta/5+\varepsilon} h=\mathrm{o}(1)$
for some $\varepsilon>0$, then equation~(\ref{cgl3.5}) remains true.
\end{thm}
\begin{rem}\label{rem3.3}
The asymptotic distribution in (\ref{cgl3.5}) is
similar to the corresponding results obtained by Karlsen \textit{et al.}
\cite{KarMykTjs07} and
Wang and Phillips \cite{WanPhi09N1}. The rate of convergence is slower
than that for the stationary time series case as $\sum_{t=1}^n
K (\frac{V_t - v}{h} )=\mathrm{O}_P(N(n)h)$ and~$N(n)$ is usually
smaller than $n$ almost surely. The condition $n^{\beta/5+\varepsilon}
h=\mathrm{o}(1)$ makes sure that the bias term of the nonparametric estimator
$\widehat{g}_n(v)$ is negligible.
\end{rem}

\subsection{Some extensions}\label{s3.2}

In this section, we give some detailed discussion of the possible
extensions raised in Remark \ref{rem2.1}(ii) and (iii). In addition, we also
suggest some other extensions.

Instead of considering a variety of extensions of model (\ref{cgl1.3})
and Theorems~\ref{thm3.1} and~\ref{thm3.2},
this section considers several extensions
that are naturally based on the relaxation of A1--A3 to
Assumptions \ref{as3.1}--\ref{as3.3} below, respectively. As a consequence, the
extended models proposed below allow for the inclusion of endogeneity,
heteroscedasticity and deterministic trending.
\begin{ass}\label{as3.1}
There are a bivariate function $H(\cdot, \cdot)$
and a stationary process $\{U_t\}$ such that
$X_t = H (V_t, \frac{t}{n} ) + U_t$ for $1\leq t\leq n$.
\end{ass}
\begin{ass}\label{as3.2}
\textup{(i)} Let \textup{A2(i)} hold.

\textup{(ii)} Let $\{\epsilon_t\}$ be of the form of either $\epsilon_t = \sigma
(\zeta_t) e_t$ or $\epsilon_t = \lambda(\xi_t) + e_t$ with $\zeta_t =
U_t$ or $V_t$ and $\xi_t = U_t$ or $\xi_t=\frac{t}{n}$, in which
$\{e_t\}$ is a stationary ergodic Markov process satisfying~\textup{A2(ii)}
and both $\sigma(\cdot)$ and $\lambda(\cdot)$ are smooth functions.
\end{ass}
\begin{ass}\label{as3.3}
\textup{(i)} Let \textup{A3(i)} hold.

\textup{(ii)} Let $\{V_t\}$ be independent of both $\{U_t\}$ and $\{e_t\}$. In
addition, $E[e_t|U_t]=0$.
\end{ass}

While it is difficult to consider some general non-stationarity for
$\{X_t\}$, it is possible to consider a general inhomogeneous case in
Assumption \ref{as3.1} to allow for a bivariate functional form of
$H(\cdot,\cdot)$ such that the non-stationarity of $\{X_t\}$ is caused by both
the involvement of $\{V_t\}$ and the dependence on $t$. In this case,
$H(\cdot, \cdot)$ may be estimated nonparametrically by
%
\begin{equation}\label{eq3.1}
\widehat{H}(v, \tau) = \sum_{t=1}^n W_{nt}(v, \tau) X_t
\qquad   \mbox{with }
W_{nt}(v, \tau) = \frac{K_{v, \tau}(V_t, t)}{\sum_{k=1}^nK_{v, \tau}(V_k, k)},
\end{equation}
where $K_{v, \tau} (V_t, t) = \frac{1}{h_1} \frac{1}{h_2} K_1 (\frac
{V_t - v}{h_1} ) K_2 (\frac{{t}/{n} - \tau}{h_2} )$,
in which both $K_i(\cdot)$ are probability kernel functions and $h_i$
are bandwidth parameters for $i=1,2$.

Assumption \ref{as3.2}(ii) allows for inclusion of endogeneity,
heteroscedasticity and deterministic trending. In the case where we
have either $\epsilon_t = \sigma(U_t)  e_t$ or
$\epsilon_t = \sigma(V_t)  e_t$ with $E[e_t|U_t]=E[e_t|V_t]=0$,
it follows that either
$E[\epsilon_t|V_t] = E[\sigma(V_t) e_t|V_t] =\sigma(V_t) E[e_t|V_t] = 0
= E[\epsilon_t]$ or $E[\epsilon_t|V_t] = E[\sigma(U_t) e_t] = E[\epsilon_t]$.
This implies Assumption \ref{as3.2}(ii) holds in both cases. In addition,
Assumption \ref{as3.2}(ii) also includes the case where $\epsilon_t = \lambda
 (\frac{t}{n} ) + e_t$ or $\epsilon_t = \lambda(U_t) + e_t$.
In such cases, obviously we have $E[\epsilon_t|V_t] = E[\epsilon_t]$.

Under Assumptions \ref{as3.1}--\ref{as3.3},
model (\ref{cgl1.3}) can be written as either
%
\begin{eqnarray}\label{eq3.2a}
Y_t & = & X_t^{\tau} \theta_0 + g(V_t) + \sigma(\zeta_t) e_t,
\nonumber\\[-8pt]
\\[-8pt]
X_t & = & H \biggl(V_t, \frac{t}{n} \biggr) + U_t,
\nonumber
\end{eqnarray}
where $\zeta_t = U_t$ or $V_t$, or
%
\begin{eqnarray}\label{eq3.2b}
Y_t & = & X_t^{\tau} \theta_0 + g(V_t) + \lambda(\xi_t) + e_t,
\nonumber\\[-8pt]
\\[-8pt]
X_t & = & H \biggl(V_t, \frac{t}{n} \biggr) + U_t,
\nonumber
\end{eqnarray}
where $\xi_t= U_t$ or $\xi_t=\frac{t}{n}$.

Estimation of $\theta_0$ and $g(\cdot)$ in (\ref{eq3.2a}) is similar to
what has been proposed in Section~\ref{s2}. Since model (\ref{eq3.2b}) is a
semi-parametric additive model, one will need to estimate $\theta_0$
based on the form $Y_t = X_t^{\tau} \theta_0 + G(V_t, \xi_t) + e_t$
with $G(v, \tau) = g(v) + \lambda(\tau)$ before both~$g(\cdot)$
and~$\lambda(\cdot)$ can be individually estimated using the marginal
integration method as developed in Section 2.3 of Gao~\cite{Gao07}.

In both cases, one will need to replace $\{w_{nt}(v)\}$ in
(\ref{eq2.1}) and $p_n(v)$ in (\ref{cgl2.12}) by $\{W_{nt}(v, \tau)\}$ of
(\ref{eq3.1}) and $p_n(v, \tau) = \frac{1}{N(n)} \sum_{k=1}^n K_{v, \tau
}(V_k, k)$, respectively.

Since the establishment and the proofs of the corresponding results of
Theorems~\ref{thm3.1} and~\ref{thm3.2}
for models (\ref{eq3.2a}) and (\ref{eq3.2b})
involve more technicalities than those given in Appendices~B and C of
the supplemental document, we wish to leave the discussion of
models~(\ref{eq3.2a}) and~(\ref{eq3.2b}) to a future paper.

\section{Simulation study}\label{s4}

To illustrate our estimation procedure, we consider a simulated example
and a real data example in this section. Throughout the section, the
uniform kernel $K(v) = \frac{1}{2} I_{[-1,1]}(v)$ is used. A difficult
problem in simulation is the choice of a proper bandwidth. From the
asymptotic results in Section~\ref{s3}, we can find that the rates of convergence
are different from those in the stationary case with $n$ being replaced
by $N(n)$. In practice, we have found it useful to use a
semi-parametric cross-validation method (see, e.g., Section
2.1.3 of H\"{a}rdle \textit{et al.} \cite{HarLiaGao00}).
\begin{exm}\label{ex4.1}
Consider a partially linear time series model of the form
%
\begin{equation}\label{cgl4.1}
Y_t=X_t\theta+ g(V_t) + \epsilon_t, \qquad  t=1,2,\ldots,n,
\end{equation}
where $V_t=V_{t-1} + v_t$ with $V_0=0$ and $\{v_t\}$ is a sequence of
i.i.d. random variables generated from $N(0,0.1^2)$, $\{\epsilon_t\}$
is generated by an AR$(1)$ model of the form
\[
\epsilon_t=0.5\epsilon_{t-1}+\eta_t,
\]
in which $\{\eta_t\}$ is a sequence of i.i.d. random variables
generated from $N(0,1)$, $\{v_t\}$ and~$\{\eta_t\}$ are mutually
independent. We then choose the true value of $\theta$ as $\theta_0=1$,
the true form of $g(\cdot)$ as $g_0(v)=v$ and consider the following
cases for $\{X_t\}$.
\begin{longlist}[(ii)]
\item[(i)] $X_t=U_t$, where $\{U_t\}$ is a sequence of i.i.d. $N(0,1)$
random variables,
\item[(ii)] $X_t=V_t+U_t$, where $\{U_t\}$ is defined as in case (i).
\end{longlist}

\begin{table}
\caption{Simulation results for the estimator of $\theta_0$\label{t1}}
\begin{tabular*}{\textwidth}{@{\extracolsep{\fill}}d{4.0}lll@{}}
\hline
\multicolumn{1}{@{}l}{$n$} & $H(\cdot)$ & AE & {SE}\\
\hline
200 &$H(v)\equiv0$& 0.0137 & 0.0144 \\
700 &$H(v)\equiv0$& 0.0117 & 0.0086 \\
1200 &$H(v)\equiv0$& 0.0064 & 0.0062 \\
200 &$H(v)=v$ & 0.0172 & 0.0215 \\
700 &$H(v)=v$ & 0.0149 & 0.0126\\
1200 &$H(v)=v$ & 0.0079 & 0.0108\\
\hline
\end{tabular*}
\end{table}
\begin{table}[b]
\caption{Simulation results for the estimator of $g_0(v)=v$\label{t2}}
\begin{tabular*}{\textwidth}{@{\extracolsep{\fill}}d{4.0}lll@{}}
\hline
\multicolumn{1}{@{}l}{$n$} & $H(\cdot)$ & AE & {SE}\\
\hline
200 & $H(v)\equiv0$ & 0.1158 & 0.0575 \\
700 & $H(v)\equiv0$ & 0.0894 & 0.0341 \\
1200 & $H(v)\equiv0$ & 0.0628 & 0.0210 \\
200 & $H(v)=v$ & 0.1391 & 0.0582 \\
700 & $H(v)=v$ & 0.1299 & 0.0437 \\
1200 & $H(v)=v$ & 0.1075 & 0.0367 \\
\hline
\end{tabular*}
\end{table}

\begin{figure}[t]

\includegraphics{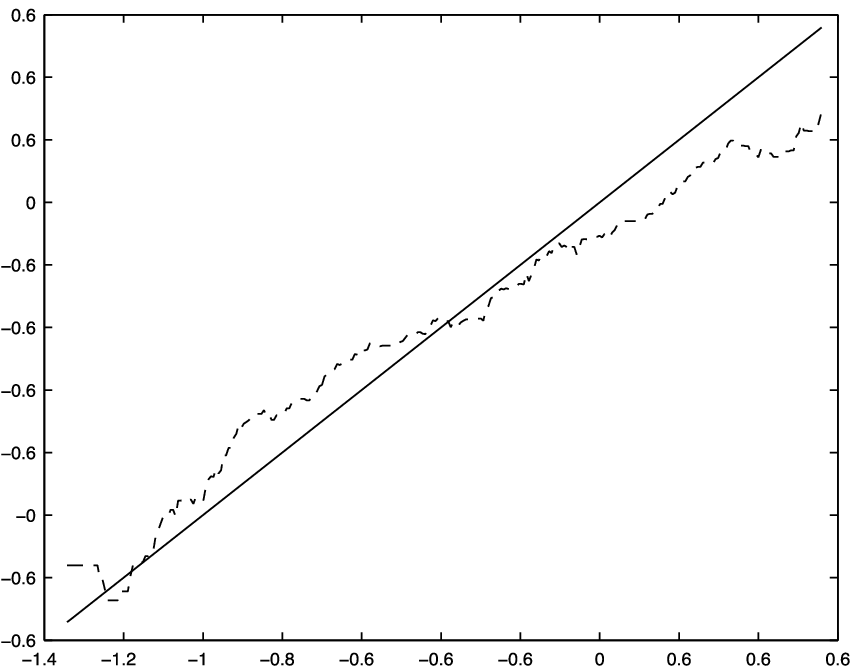}

 \caption{Nonparametric estimate of the regression function $g_0(v)$ for the
case of $H(v) \equiv 0$~with sample size $n = 200$; the solid line is the
true line, and the dashed curve is the estimated curve.}\label{f1}
\end{figure}
\begin{figure}[b]

\includegraphics{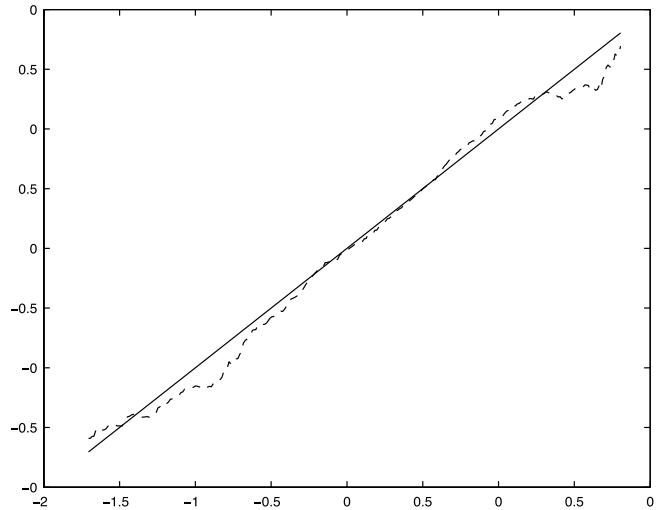}

 \caption{Nonparametric estimate of the regression function $g_0(v)$ for
the case of $H(v) \equiv 0$~with sample size $n = 700$; the solid line is the
true line, and the dashed curve is the estimated curve.}\label{f2}
\end{figure}

It is easy to check that the random walk $\{V_t\}$ defined in this
example corresponds to a $1/2$ null recurrent process and the
assumptions in Section \ref{s2} are satisfied here. We choose sample
sizes $n=200,  700, 1200$ and $N=1000$ as the number of
replications in the simulation. The simulation results are listed
in Tables \ref{t1} and \ref{t2} and the plots are given in Figures \ref{f1}--\ref{f6}.

\begin{figure}

\includegraphics{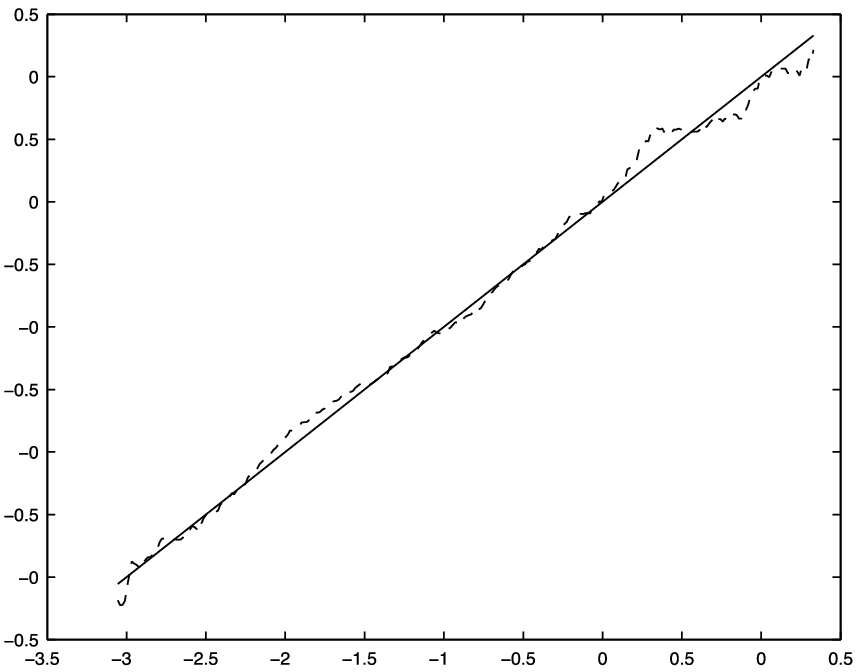}

 \caption{Nonparametric estimate of the regression function $g_0(v)$ for
the case of $H(v) \equiv 0$~with sample size $n = 1200$; the solid line is
the true line, and the dashed curve is the estimated curve.}\label{f3}
\end{figure}
\begin{figure}[b]

\includegraphics{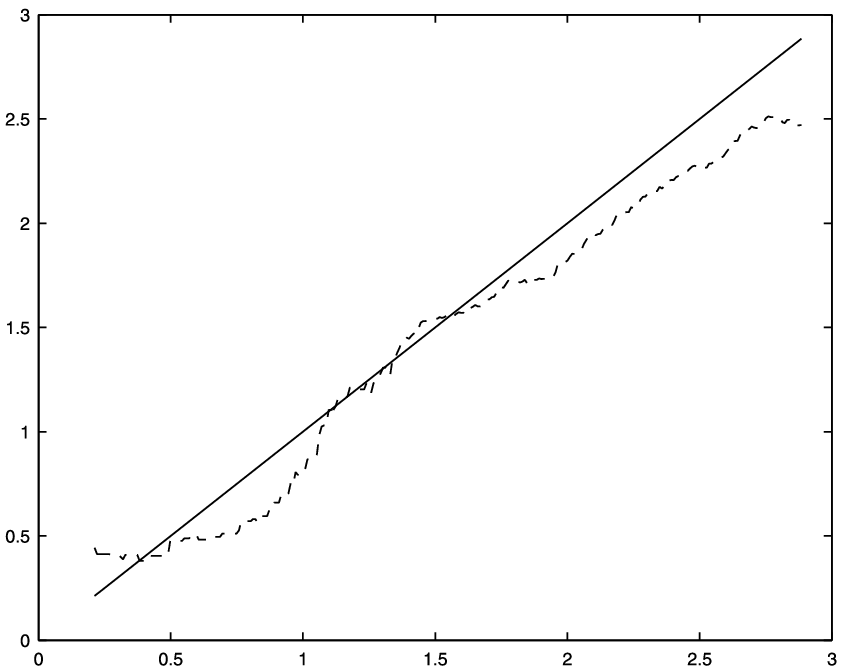}

 \caption{Nonparametric estimate of the regression function $g_0(v)$ for
the case of $H(v) = v$~with sample size $n = 200$; the solid line is the
true line, and the dashed curve is the estimated curve.}\label{f4}
\end{figure}
\begin{figure}

\includegraphics{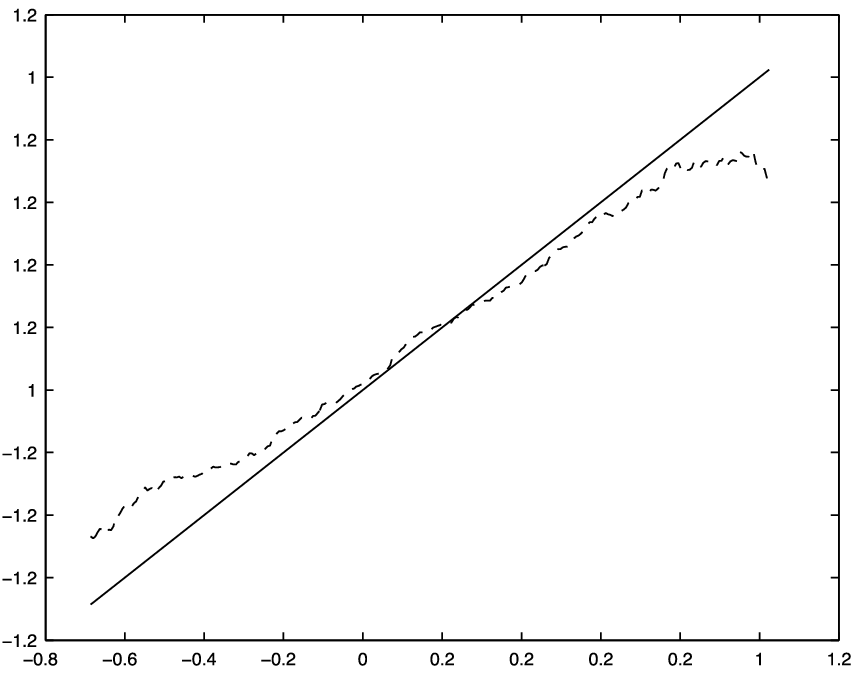}%
\vspace*{-3pt}
 \caption{Nonparametric estimate of the regression function $g_0(v)$ for
the case of $H(v) = v$~with sample size $n = 700$; the solid line is the
true line, and the dashed curve is the estimated
curve.}\label{f5}\vspace*{-3pt}
\end{figure}
\begin{figure}[b]
\vspace*{-3pt}
\includegraphics{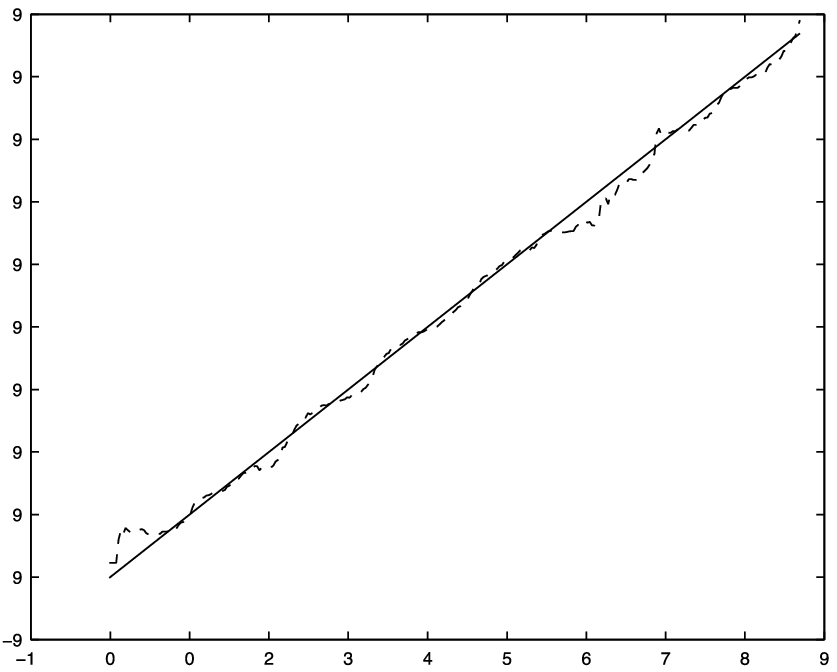}%
\vspace*{-3pt}
 \caption{Nonparametric estimate of the regression function $g_0(v)$ for
the case of $H(v) = v$~with sample size $n = 1200$; the solid line is
the true line, and the dashed curve is the estimated curve.}\label{f6}
\end{figure}

The performance of $\widehat{\theta}_n$ is given in Table \ref{t1}.
The ``AE'' in Table \ref{t1} is\vspace*{2pt} defined by
$\frac{1}{1000}\sum_{j=1}^{1000}  |\widehat{\theta}(j)-\theta_0 |$,
where $\widehat{\theta}(j)$ is the value of $\widehat{\theta}_n$
in the $j$th replication. ``SE'' is the standard error of
$\{\widehat{\theta}(j)\}$. From Table \ref{t1}, we find that the
estimator of $\theta_0$ performs well in the small and
medium sample cases and it improves when the sample size
increases.\looseness=1

The performance of the nonparametric estimator is given in Table \ref{t2}. The
``AE'' in Table~\ref{t2} is the mean of the
absolute errors in 1000 replications. The absolute error is
defined by $\frac{1}{300}\sum_{j=1}^{300} |\widehat{g}_n(v_j)-v_j |$,
where $v_j=v_{\min}+\frac{j-1}{300}(v_{\max}-v_{\min})$ for
$j=1,2,\ldots,300$,~$v_{\max}$ and $v_{\min}$ are the maximum and
minimum of the random walk $\{V_t,1\leq t\leq n\}$, respectively.
``SE'' in Table \ref{t2} is the standard error. From Table \ref{t2}, we find that the
nonparametric estimate of $g_0(v)=v$
performs well in our example and it improves when the sample size
increases.\vadjust{\goodbreak}

Figures \ref{f1}--\ref{f3} compare the true nonparametric regression function
$g_0(\cdot)$ and its nonparametric estimator for the case of
$H(v)=0$ when the sample sizes are 200, 700 and 1200, respectively.
Figures \ref{f4}--\ref{f6} compare the true nonparametric regression function
$g_0(\cdot)$ with its nonparametric estimator for the case of
$H(v)=v$ when the sample sizes are 200, 700 and 1200, respectively.
The solid line is $g_0(\cdot)$ and the dashed line is the
nonparametric estimator. We cannot forecast the trace of the
random walk $\{V_t\}$ because of its non-stationarity. Hence, we
estimate the true regression function $g_0(\cdot)$ according to
the scope of $\{V_t\}$ and we cannot estimate $g_0(\cdot)$ in
other points out of the scope since there is not enough sample in
the neighborhood of each of such points. That is why the scopes of the abscissa
axis are different in Figures \ref{f1}--\ref{f6}. We can also find that the
performance of the nonparametric estimate of $g_0(\cdot)$
improves as the sample size increases.\looseness=-1
\end{exm}

\section{An empirical application}\label{s5}

\begin{figure}[b]

\includegraphics{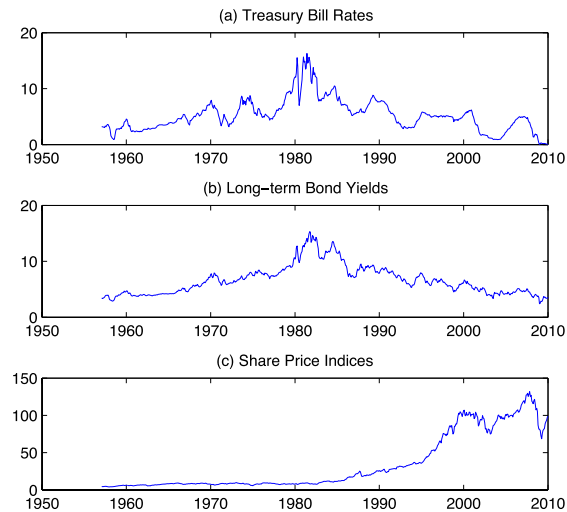}

 \caption{Time plots of the three series used in
Section \protect\ref{s5} over the period of Jan/1957--Dec/2009 with 624~observations.
(a) treasury bill rates; \textup{(b)} long-term bond yields; (c) share
prices.\label{f7}}
\end{figure}

We use monthly observations on the U.S. share price indices, long-term
government bond yields and treasury bill rates from Jan/1957--Dec/2009.
The data are obtained from the International Monetary Fund's (IMF)
International Financial Statistics (IFS). The share price series used
is IFS Series 11162ZF. The long-term government bond yield, which is
the 10-year yield, is from the IFS Series 11161ZF. The treasury bill
rate is from IFS Series 11160CZF. Figure \ref{f7}(a)--(c) gives the data plots
of the share prices, the long-term bond yields and the treasury bill rates.

To see whether there exist some statistical evidences for the three
series to have the unit root type of non-stationarity, we carry out a
Dickey--Fuller (DF) unit root test on the three series. We first fit
the data by an $\operatorname{AR}(1)$ model of the form
\[
Z_{t}=\rho Z_{t-1}+e_{t},
\]
where $Z_t= $ share price at time $t$ or long-term bond yield at time
$t$ or treasury bill rate at time~$t$.
Then, by using the least-squares estimation method, we estimate the
parameter~$\rho$ for the three series: for the share price series,
$\widehat{\rho}_{\rm share}=1.0023$; for the long-term bond yield
series, $\widehat{\rho}_{\mathrm{Lbond}}=0.9992$; and for the treasury bill
rate series, $\widehat{\rho}_{\mathrm{Tbill}}=0.9966$. Then we calculate the
Dickey--Fuller $t$ statistics and compare them with the critical values
at the $5\%$ significance level. The simulated $P$ values for the
long-term bond yields, treasury bill rates and share prices are
$0.7040$, $0.3130$ and $0.4410$, respectively. In addition, we also
employ an augmented DF test and the nonparametric test proposed in Gao
\textit{et al.}~\cite{Gaoetal09} for checking the unit root structure of $\{Z_t\}$.
The resulting $P$ values are very similar to those obtained above.

\begin{figure}

\includegraphics{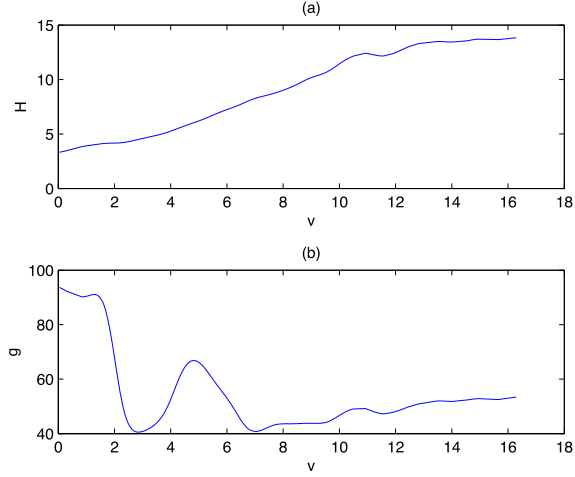}%
\vspace*{-3pt}
\caption{Estimates of the nonparametric functions
$H(\cdot)$ and $g(\cdot)$ in Case A.\label{f8}}
\vspace*{-3pt}
\end{figure}

Therefore, both the estimation results and the simulated $P$ values
suggest that there is some strong evidence for accepting the null
hypothesis that a unit root structure exists in these series at the $5\%$
significance level.

We then consider the following modelling problem:
\begin{eqnarray*}
Y_t & = & X_t \theta_0 + g(V_t) + \epsilon_t,
\\[-2pt]
X_t & = & H(V_t) + U_t,
\end{eqnarray*}
where Case A: $Y_t$ is the share price, $X_t$ is the long-term bond
yield and $V_t$ is the treasury bill; and Case B: $Y_t$ is the
long-term bond yield, $X_t$ is the share price and $V_t$ is the
treasury bill.

\begin{figure}

\includegraphics{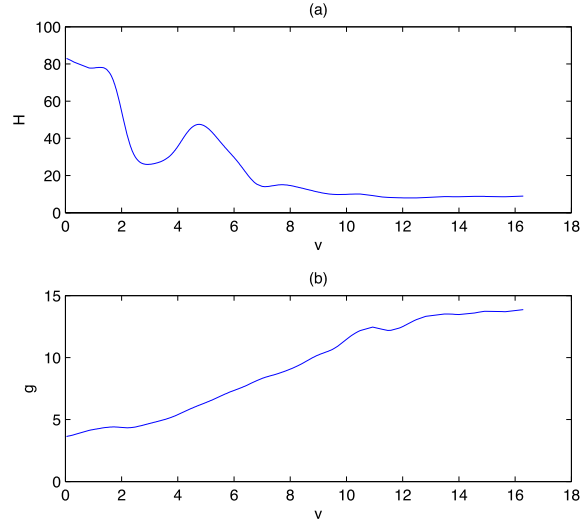}

 \caption{Estimates of the nonparametric functions
$H(\cdot)$ and $g(\cdot)$ in Case B.\label{f9}}
\end{figure}

For Case A, the resulting estimator of $\theta_0$ is
$\widehat{\theta}=-3.2155$ and the plots of the estimates of $g(\cdot)$ and $H(\cdot)$
are given in Figure \ref{f8}. For Case B, the resulting\vadjust{\goodbreak} estimator of $\theta
_0$ is $\widehat{\theta}=-0.0037$ and the the plots of the estimates of
$g(\cdot)$ and $H(\cdot)$ are given in Figure~\ref{f9}.\looseness=1

Figures \ref{f8} and \ref{f9} show that increases in treasury bill rates tend to lead
to increases in long-term bond yields and decreases in share prices.
Such findings are supported by the theory of finance and consistent
with existing studies. Moreover, Figures \ref{f7}--\ref{f9} clearly indicate our new
findings that both null recurrent non-stationarity and nonlinearity can
be simultaneously exhibited in the share price, the long-term bond
yield and the treasury bill rate variables.

Due to the cointegrating relationship among the stock price, the
treasury bill rate and the long-term bond yield variables, our
experience suggests that models (\ref{eq3.2a}) and (\ref{eq3.2b}) might
be more suitable for this empirical study. We will have another look at
this data after models (\ref{eq3.2a}) and~(\ref{eq3.2b}) have been
fully studied.

\section{An outline of the proofs of the theorems}\label{s6}

In this section, we provide only one key lemma and then an outline of
the proofs of Theorems \ref{thm3.1} and \ref{thm3.2}. The detailed proofs of the theorems
are available from the supplemental document by
Chen, Gao and Li \cite{CheGaoLi10}.
\begin{lem}\label{lem6.1}
Under the conditions of Theorem \ref{thm3.1},
we have as $n\rightarrow\infty$,
%
\begin{equation}\label{cglB.2}
\frac{1}{n}\widetilde{X}^\tau\widetilde{X}
\stackrel{P}{\longrightarrow}\Sigma.
\end{equation}
\end{lem}
\begin{pf*}{Proof of Theorem \ref{thm3.1}}
In view of Lemma \ref{lem6.1} and the decomposition
\begin{eqnarray*}
\widetilde{X}^{\tau}\widetilde{X} (\widehat{\theta}_n-\theta_0 )
&=&\widetilde{X}^\tau(\widetilde{Y}-\widetilde{X}\theta_0)
\\[-2pt]
&=&\sum_{t=1}^{n}\widetilde{X}_{t}\widetilde{g}(V_{t})F_t
+ \sum_{t=1}^{n}\widetilde{X}_{t}\epsilon_{t}F_t
- \sum_{t=1}^{n}\widetilde{X}_{t}F_t
\Biggl(\sum_{k=1}^{n}w_{nk}(V_{t})\epsilon_{k} \Biggr),
\end{eqnarray*}
in order to prove Theorem \ref{thm3.1}, we need only to show that for large
enough $n$
%
\begin{eqnarray}
\sum_{t=1}^{n}\widetilde{X}_{t}\widetilde{g}(V_{t})F_t
& = & \mathrm{o}_P\bigl(\sqrt{n}\bigr),
\label{eq:jiti1}\\[-2pt]
\sum_{t=1}^{n}\widetilde{X}_{t}F_t
\Biggl\{\sum_{k=1}^{n}w_{nk}(V_{t}) \epsilon_{k} \Biggr\}
& = & \mathrm{o}_P\bigl(\sqrt{n}\bigr),
\label{eq:jiti2}\\[-2pt]
n^{-1/2}\sum_{t=1}^{n}\widetilde{X}_{t}\epsilon_{t}F_t
& \stackrel{d}{\longrightarrow} &
N (0, \Sigma_{\epsilon,U} ),
\label{eq:jiti3}
\end{eqnarray}
where $\widetilde
{g}(V_t)=g(V_t)-\sum_{k=1}^nw_{nk}(V_t)g(V_k)$.
Recall that\vspace*{2pt} $\widetilde{X}_t = X_t - \sum_{s=1}^n w_{ns}(V_t) X_s = U_t
- \sum_{s=1}^n w_{ns}(V_t) U_s + \widetilde{H}(V_t)$, where
$\widetilde{H}(V_t)=H(V_t) - \sum_{s=1}^nw_{ns}(V_t) H(V_s)$.

In order to prove (\ref{eq:jiti1})--(\ref{eq:jiti3}), it suffices to
show that for large enough $n$
%
\begin{eqnarray}
\sum_{t=1}^{n}U_t \widetilde{g}(V_{t})F_t & = & \mathrm{o}_P\bigl(\sqrt{n}\bigr),
\label{eq:jiti4a}\\[-2pt]
\sum_{t=1}^{n} \overline{U}_t\widetilde{g}(V_{t}) F_t
& = & \mathrm{o}_P\bigl(\sqrt{n}\bigr),
\label{eq:jiti4b}\\[-2pt]
\sum_{t=1}^{n} \widetilde{g}(V_{t})\widetilde{H}(V_{t}) F_t
& = & \mathrm{o}_P\bigl(\sqrt{n}\bigr),
\label{eq:jiti4}\\[-2pt]
\sum_{t=1}^{n} U_t\overline{\epsilon}_t F_t
& = & \mathrm{o}_P\bigl(\sqrt{n}\bigr),
\label{eq:jiti5a}\\[-2pt]
\sum_{t=1}^{n} \overline{U}_t\overline{\epsilon}_t F_t
& = & \mathrm{o}_P\bigl(\sqrt{n}\bigr),
\label{eq:jiti5}\\[-2pt]
\sum_{t=1}^{n} \widetilde{H}(V_{t})\overline{\epsilon}_t F_t
& = & \mathrm{o}_P\bigl(\sqrt{n}\bigr),
\label{eq:jiti6}\\[-2pt]
\sum_{t=1}^{n} \overline{U}_t \epsilon_{t} F_t
& = & \mathrm{o}_P\bigl(\sqrt{n}\bigr),
\label{eq:jiti6a}\\[-2pt]
\sum_{t=1}^{n} \widetilde{H}(V_{t}) \epsilon_t F_t
& = & \mathrm{o}_P\bigl(\sqrt{n}\bigr),
\label{eq:jiti7}\\[-2pt]
n^{-1/2}\sum_{t=1}^{n}U_t
\epsilon_{t}F_t&\stackrel{d}{\longrightarrow}&
N (0, \Sigma_{\epsilon,U} ),\vspace*{-2pt}
\label{eq:jiti8}
\end{eqnarray}
where
$\overline{U}_t = \sum_{s=1}^n w_{ns}(V_t) U_s$ and
$\overline{\epsilon}_t = \sum_{s=1}^n w_{ns}(V_t)
\epsilon_s$.

In the following, we verify equations
(\ref{eq:jiti4a})--(\ref{eq:jiti8}) to complete the proofs of
Theorem~\ref{thm3.1}(i) and Theorem~\ref{thm3.1}(ii). Note that, for Theorem
\ref{thm3.1}(i), equations (\ref{eq:jiti4}), (\ref{eq:jiti6}) and~(\ref{eq:jiti7}) hold trivially.

By the continuity of $g(\cdot)$ and $g^{\prime}(\cdot)$, we have for
$1\leq t\leq n$,
%
\begin{eqnarray}\label{appendix1}
&& \frac{1}{N(n)h}\sum_{j=1}^n K \biggl(\frac{V_j-V_t}{h} \biggr)
\bigl(g(V_j)-g(V_t)\bigr)
\nonumber\\[-9.5pt]
\\[-9.5pt]
&&\quad = \frac{g^{\prime}(V_t)}{N(n)h}\sum_{j=1}^n(V_j-V_t) K
\biggl(\frac{V_j-V_t}{h} \biggr) \bigl(1+\mathrm{o}_P(1)\bigr).
\nonumber\vspace*{-2pt}
\end{eqnarray}

Thus, in view of (\ref{appendix1}) and Lemma 3.4 of
Karlsen and Tj\o stheim \cite{KarTjs01},
in order to prove~(\ref{eq:jiti4a}), it suffices to
show that for $n$ large enough
%
\begin{equation}\label{appendix2}
\sum_{t=1}^{n}U_t \Delta_n(V_t) F_t = \mathrm{o}_P\bigl(\sqrt{n}\bigr),\vspace*{-2pt}
\end{equation}
where $\Delta_n(V_t) = \frac{g^{\prime}(V_t)}{n^{\beta- \eta}h
p_n(V_t)} \sum_{j=1}^n  (V_j - V_t ) K (\frac{V_j -V_t}{h} )$.

This kind of procedure of replacing $N(n)$ by $n^{\beta-\eta}$ and
ignoring a small-order term as involved in (\ref{appendix1}) will be
used repeatedly throughout the proofs in Appendices B and C of the
supplemental document.

We then may show that (\ref{eq:jiti4b}) holds. Similarly to
(\ref{appendix1}) and (\ref{appendix2}), we need only to show that
%
\begin{equation}\label{appendix3}
\sum_{t=1}^{n}\widehat{U}_t\Delta_n(V_t) F_t=\mathrm{o}_P\bigl(\sqrt{n}\bigr),\vspace*{-2pt}
\end{equation}
where $\widehat{U}_t=\frac{1}{n^{\beta-\eta}h  p_n(V_t)}
(\sum_{k=1}^nK (\frac{V_k-V_t}{h} )U_k )$.\vspace*{2pt}

The detailed derivations for (\ref{appendix2}) and (\ref{appendix3})
are available from Appendix B of the supplemental document. The
detailed proofs of (\ref{eq:jiti5a}), (\ref{eq:jiti5}),
(\ref{eq:jiti6a}) and (\ref{eq:jiti8}) are also available from Appendix B.
This will complete the proof of Theorem \ref{thm3.1}(i).

We then may prove Theorem \ref{thm3.1}(ii) by completing the proofs of
(\ref{eq:jiti4}), (\ref{eq:jiti6}) and~(\ref{eq:jiti7}), which are again
available from Appendix B of the supplemental document.\vspace*{-2pt}
\end{pf*}
\begin{pf*}{Proof of Theorem \ref{thm3.2}}
By the definition of $\widehat{g}_n(v)$, we have
%
\begin{eqnarray}\label{cglB.43}
\widehat{g}_n(v)-g(v)&=&\sum_{t=1}^{n}
w_{nt}(v)(Y_{t}-X_{t}\widehat{\theta}_n)-g(v)
\nonumber\\[-9.5pt]
\\[-9.5pt]
& = & \sum_{t=1}^{n} w_{nt}(v)\bigl(\epsilon_t + g(V_t)-g(v)\bigr)
+ \sum_{t=1}^{n} w_{nt}(v)X_{t}(\theta_0-\widehat{\theta}_n).
\nonumber\vspace*{-2pt}\vadjust{\goodbreak}
\end{eqnarray}

Let $\Phi_{n,1} = \sum_{t=1}^{n} w_{nt}(v)(\epsilon_t +
g(V_t)-g(v)$ and $\Phi_{n,2}=\sum_{t=1}^{n} w_{nt}(v)X_{t}(\theta
_0-\widehat{\theta}_n)$. Then, we have
%
\begin{eqnarray}\label{cglB.44}
\widehat{g}_n(v)-g(v)&=&\sum_{t=1}^{n}
w_{nt}(v)(Y_{t}-X_{t}\widehat{\theta}_n)-g(v)
= \Phi_{n,1} + \Phi_{n,2}.
\end{eqnarray}

Since $\{\epsilon_t\}$ is assumed to be stationary and $\alpha$
mixing, by Corollary 5.1 of Hall and Heyde~\cite{HalHey80} and an existing
technique to deal with the bias term (see, e.g., the proof of
Theorem 3.5 of Karlsen \textit{et al.} \cite{KarMykTjs07}),
we have as $n\rightarrow\infty$
%
\begin{equation}\label{cglB.45}
\sqrt{\sum_{t=1}^n K \biggl(\frac{V_t -v}{h} \biggr)}
\Phi_{n,1}\stackrel{d}{\longrightarrow}N
\biggl(0,\sigma^{2}\int K^2(u)\,\mathrm{d}u \biggr).
\end{equation}

By (\ref{cglB.43})--(\ref{cglB.45}), it is sufficient to show that
%
\begin{equation}\label{cglB.46}
\sqrt{\sum_{t=1}^n K \biggl(\frac{V_t-v}{h} \biggr)}
\Phi_{n,2}=\mathrm{o}_P(1).
\end{equation}

The proof of (\ref{cglB.46}) may then be completed by Theorem \ref{thm3.1} and
Assumptions A1--A6. The details are available from Appendix B of the
supplemental document. This completes an outline of the proofs of
Theorems \ref{thm3.1} and \ref{thm3.2}.
\end{pf*}

\section*{Acknowledgements}

This work was started when the first and third authors were visiting
the second author in 2006/2007. The authors would all like to thank the
Editor, the Associate Editor and two references for their constructive
comments on an earlier version. Thanks also go to the Australian
Research Council Discovery Grants Program for its financial support
under Grant Numbers DP0558602 and DP0879088.

\begin{supplement}
\stitle{Proofs of the theorems}
\slink[doi]{10.3150/10-BEJ344SUPP} 
\sdatatype{.pdf}
\sfilename{BEJ344\_SUPPL.pdf}
\sdescription{We provide this supplemental document in case the reader may want to
have a look at the detailed proofs of Theorems \ref{thm3.1}~and~\ref{thm3.2}
and Lemma \ref{lem6.1}.
The details are available from Chen, Gao and Li \cite{CheGaoLi10}.}
\end{supplement}

%

\printhistory

\end{document}